\documentclass{scrartcl}

\usepackage{import}
\usepackage{todonotes}
\usepackage{lmodern}
\usepackage[utf8]{inputenc}
\usepackage[margin=1.1in]{geometry}
\usepackage{amsmath,amsthm,amssymb,amsfonts,color,dsfont,enumitem,soul}
\usepackage[hidelinks]{hyperref}
\usepackage[normalem]{ulem}
\usepackage{tikz,tikz-cd}
\usepackage{multicol}
\usepackage{float}
\usepackage{siunitx,booktabs,array}
\usepackage{epigraph}
\usepackage{pdflscape}
\usepackage[capitalise]{cleveref}

\usepackage{listings}
\usepackage{color}

\usepackage[backend=biber,isbn=false,url=false, sorting=none, style=numeric-comp]{biblatex}
\newcommand{\gr}[1]{\textcolor{gray}{#1}}

\usepackage{tcolorbox}

\tcbuselibrary{theorems}

\newtcbtheorem[number within=section]{remark}{Remark}%
{colback=black!10, colframe=black,
  coltitle=white,
  boxrule=0.8pt,
  arc=4pt, 
 }
{th}

\begin{document}
\title{AnyonWiki.org}
\subtitle{An introduction to the open-source database \\ of fusion categories and fusion rings}
\author{Colleen Delaney, Ralph Razzouk, Joost Slingerland, Gert Vercleyen}
\maketitle

\definecolor{latte}{rgb}{.99,.97,.92}
\definecolor{codeblue}{rgb}{0.08, .3, 1}
\definecolor{codegray}{rgb}{0.5,0.5,0.5}
\definecolor{codered}{rgb}{.76, 0, .22}
\definecolor{codegreen}{rgb}{0,.76,.54}
\lstdefinestyle{mystyle}{
    backgroundcolor=\color{latte},
    commentstyle=\color{codeblue},
    keywordstyle=\color{codered},
    numberstyle=\tiny\color{codegreen},
    stringstyle=\color{codegray},
    breakatwhitespace=false,
    breaklines=true,
    captionpos=b,
    keepspaces=true,
    numbers=left,
    numbersep=5pt,
    showspaces=false,
    showstringspaces=false,
    showtabs=false,
    tabsize=2
}
\lstset{style=mystyle}

\begin{abstract}

The AnyonWiki is an open-source project that provides full algebraic and numerical data of fusion and unitary modular fusion categories as a resource for the scientific community. It comes with a database, the AnyonWikiDatabase, that currently includes the classification of multiplicity-free pivotal (braided) fusion categories through rank 7 and more than 25000 fusion rings including the classification of multiplicity-free fusion rings through rank 10, with work in progress to expand. This data is available through AnyonWiki.org, a website that provides robust search and download functionality and presents each fusion ring and category via a dedicated webpage. In this short article we discuss the data available at AnyonWiki.org as well as the scientific and communal motivations for the design of its underlying AnyonWikiDatabase, which are centered around making the data accessible, independently verifiable, granularly citable, useful, and easy to work with.

\end{abstract}

\section{Introduction}
\label{sec:introduction}

Fusion categories generalize the representation theory of finite groups and provide a unifying mathematical framework for classifying (1) anyons and defects in gapped topological phases of matter \cite{Kitaev2006, Levin2005,BBCW}, (2) topological quantum field theories \cite{TuraevViro, BarrettWestbury, ReshetikhinTuraev:invariants}, and (3) non-invertible symmetries \cite{Aasen:2016, Thorngren:2019, Aasen:2020}, as well as understanding each of their applications to quantum computation and error correction \cite{Rowell2018, KKR2010}. 

Unlike more general kinds of tensor categories, fusion categories admit a finite description as the solutions to a finite set of algebraic equations, and can therefore be stored on a computer. Finding solutions to these equations by computer, which is particularly tractable for \emph{low-rank}, \emph{multiplicity-free} fusion categories, is an exciting frontier in mathematics and theoretical physics. 

In spite of --- and perhaps because of --- the important role they play in several different areas of math and physics, fusion categories can be challenging to work with. Widely varying jargon, differing symbolic conventions, and the natural degrees of freedom in fusion category data known as \emph{gauge equivalence} are just some of the obstacles to doing interdisciplinary and translational work using them.

AnyonWiki.org is an open-source database that provides validatable fusion ring and fusion category data in a user-friendly interface for the scientific communities that depend on them. At the time of writing, the database includes all multiplicity-free fusion rings up to rank $10$ and a census of multiplicity-free fusion (including unitary, braided, and modular fusion) categories up to rank $7$, with many plans to expand going forward.

\subsection*{History of the AnyonWiki}
Work on the public-facing side of the AnyonWiki began in 2019 while G.~V.~was a PhD student at Maynooth University with J.~S., with the first iteration of the project hosted on MediaWiki starting in November 2019. The AnyonWiki was hosted on GitHub at \href{https://anyonwiki.github.io/}{anyonwiki.github.io} from June 2023 onwards, with C.D.~joining  the project in January 2024 and R.R.~in September 2025. The project was hosted at AnyonWiki.org beginning in October 2025, with the full migration from GitHub completed in August 2026. 

\subsection*{Contents of this paper}
In \cref{sec: fusion categories and fusion rings} we give an overview of the fusion category and fusion ring data which is currently available on AnyonWiki.org. The data on the website is fetched from the AnyonWikiDatabase, whose structure and design is discussed in \cref{sec: anyonwikidatabase}. In \cref{sec: citing} we explain when and how to cite use of AnyonWiki.org and the data contained therein, as well as how references to the literature and software versions and revisions are stored within individual entries of the AnyonWikiDatabase. Finally in \cref{sec: future and contributing} we discuss the future of the AnyonWiki and how individuals can contribute to the project.

\subsection*{Acknowledgments}
We thank Eddy Ardonne, Aatmun Baxi, S\'{e}bastien Palcoux, Eric Rowell, Eric Samperon, Andrew Schopieray, and Sami Zagha for their support of the project and for providing input over the years that helped shape the AnyonWiki. Ralph Razzouk was supported in part by NSF CCF 2330130.

\subsection*{AI-independent data generation and validation}
No AI/LLMs were used in the development of the software that was used to generate the fusion category or fusion ring data on AnyonWiki.org. No AI/LLMs were used in the preparation of this paper. AI/LLMs were used to aid in front-end development of the website.

\section{The fusion categories and fusion rings on AnyonWiki.org}
\label{sec: fusion categories and fusion rings}
Every pivotal (braided) fusion category and every fusion ring on the AnyonWiki is assigned both a unique AnyonWiki code, as well as a 32-digit UUID (universally unique identifier). This is necessary for internal validation as well as for archiving our data, see also \cref{sec: anyonwikidatabase} for more details. 

\subsection{Fusion categories}
\label{sec: fusion categories}
\begin{figure}[h!]
  \begin{center}
  \includegraphics[width=.7\linewidth]{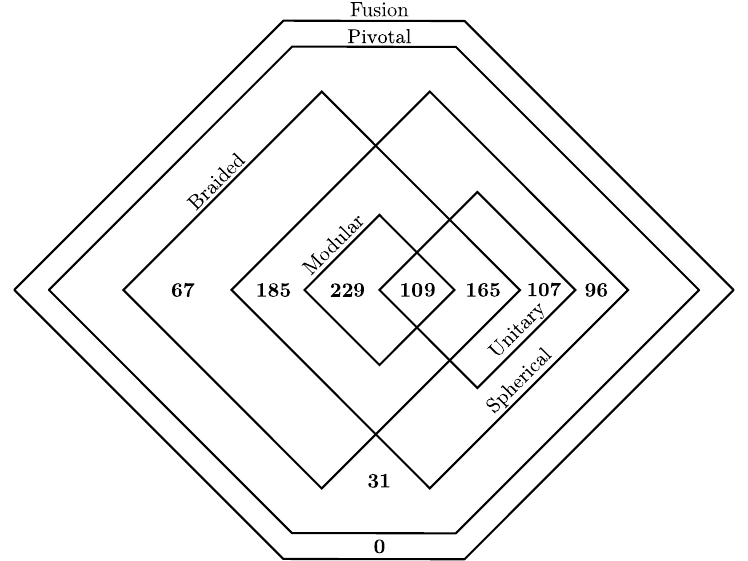}
  \end{center}
  \caption{Venn diagram of the number of multiplicity-free fusion categories up to rank $7$ with certain properties on the AnyonWiki. By definition, the ribbon fusion categories belong to the intersection of the braided and spherical squares. Each number indicates the number of fusion categories that lie in the connected region of the Venn diagram. There are thus, e.g., $0$ fusion categories with no compatible extra structures (since all the categories found admit a pivotal structure), $31$ pivotal fusion categories with no compatible extra structures, $185$ ribbon fusion categories that are neither modular nor unitary, $109$ fusion categories that are modular and unitary, $229$ fusion categories that are modular but non-unitary and thus $338$ modular fusion categories in total.}\label{fig:venn}
\end{figure}

The AnyonWiki currently contains a census of $989$ complex multiplicity-free pivotal (braided) fusion categories up to rank $7$ as obtained in~\cite{Vercleyen2026}, with a clear indication of whether each category has the property of being spherical, unitary, braidable, ribbon, and/or modular. Every available category comes with the $F$-symbols in some gauge, with the data of a specific pivotal structure (the $P$-symbols), and, if applicable, the data of a specific braided structure (the $R$-symbols). Thus two distinct entries in the database may have gauge equivalent $F$-symbols when they are inequivalent as pivotal (and possibly braided) fusion categories.
Figure~\ref{fig:venn} displays a Venn diagram of the number of categories on the AnyonWiki which admit a given structure or property.

Every pivotal (braided) fusion category has a dedicated page where visitors can find information on important properties, including:
\begin{itemize}

\item The \textbf{rank}, \textbf{multiplicity}, \textbf{Frobenius--Perron dimension}, \textbf{common names} from various sources of literature, and whether the category is \textbf{spherical, unitary, braided, ribbon, and/or modular}. The seven integer digit \textbf{AnyonWiki code} and \textbf{UUID} is also displayed.

\item The name and fusion table of the parent \textbf{fusion ring} (also known as the \textbf{Gr\"{o}thendieck ring}) of the category and a link to its individual page entry, where further data and invariants specific to fusion rings are provided as explained in the next \cref{sec: fusion rings}.

\item The \textbf{$\boldmath F$-symbols}, \textbf{pivotal coefficients} (or $P$-symbols)\footnote{Some readers may be more familiar with \textbf{Frobenius-Schur indicators}, which coincide with the pivotal coefficients of self-dual objects for unitary fusion categories, but the $P$-symbols are more general.}, and, if applicable, the \textbf{$\boldmath R$-symbols} of the pivotal (braided) fusion category, expressed in a fixed gauge both symbolically as algebraic numbers and numerically as floating point numbers.

\item The \textbf{Frobenius--Perron dimensions} of simple objects, and for spherical (including unitary but also anti-unitary) fusion categories, the \textbf{quantum dimensions} of simple objects. 

\item For modular fusion categories, the \textbf{$\boldmath S$-matrix} and \textbf{twists}, i.e.~the \emph{modular data}.

\item Whether the $\boldmath F$-, $\boldmath P$-, and $\boldmath R$-symbols are expressed in a \textbf{unitary gauge}. For any multiplicity-free category this should always be the case if the underlying fusion category is unitary.

\item A \textbf{field extension} of $\mathbb{Q}$ over which the given $F$-, $P$-, and $R$-symbols are defined, which is \emph{a priori} gauge dependent. At the moment this entry is listed as the field of algebraic numbers; in the future we will also provide the symbols defined over a gauge-invariant \emph{minimal number field} (see, e.g.,~\cite{maurer2026f} for the definition of a minimal field and how to find it). 

\end{itemize}

Additional information and new fields can easily be added over time as new properties and invariants of fusion categories are developed.

\subsection{Fusion rings}
\label{sec: fusion rings}
The AnyonWiki currently contains $25,332$ fusion rings. This includes all multiplicity-free fusion rings up to rank $10$, as well as a variety of fusion rings with higher multiplicity, all found using the methods from~\cite{v2023}. Table~\ref{tab:numpermultrank} displays the available number of fusion rings per rank and multiplicity.
\begin{table}[h!]
  \centering
\begin{tabular}{|l|rrrrrrrrrr|}
\hline  
   & 1 & 2 & 3  & 4   & 5    & 6    & 7    & 8   & 9   & 10 \\
\hline 
1  & \textbf{1} & \textbf{2} & \textbf{4}  & \textbf{10}  & \textbf{16}   & \textbf{39}   & \textbf{43}   & \textbf{96}  & \textbf{142} & \textbf{194} \\
2  & \textbf{0} & \textbf{1} & \textbf{3}  & \textbf{17}  & \textbf{37}   & \textbf{154}  & \textbf{319}  & \gr{874} &   &  \\
3  & \textbf{0} & \textbf{1} & \textbf{4}  & \textbf{24}  & \textbf{82}   & \textbf{384}  & \gr{562}  &    &    &  \\
4  & \textbf{0} & \textbf{1} & \textbf{6}  & \textbf{45}  & \textbf{134}  & \textbf{872}  & \gr{1236} &    &    &  \\
5  & \textbf{0} & \textbf{1} & \textbf{5}  & \textbf{55}  & \textbf{209}  & \gr{533}  &     &    &    &  \\
6  & \textbf{0} & \textbf{1} & \textbf{9}  & \textbf{81}  & \textbf{336}  & \gr{872}  &     &    &    &  \\
7  & \textbf{0} & \textbf{1} & \textbf{6}  & \textbf{92}  & \textbf{477}  & \gr{976}  &     &    &    &  \\
8  & \textbf{0} & \textbf{1} & \textbf{10} & \textbf{137} & \textbf{733}  & \gr{1672} &     &    &    &  \\
9  & \textbf{0} & \textbf{1} & \textbf{12} & \textbf{151} & \textbf{863}  &     &     &    &    &  \\
10 & \textbf{0} & \textbf{1} & \textbf{9}  & \textbf{186} & \textbf{1082} &     &     &    &    &  \\
11 & \textbf{0} & \textbf{1} & \textbf{10} & \textbf{238} & \textbf{1383} &     &     &    &    &  \\
12 & \textbf{0} & \textbf{1} & \textbf{20} & \textbf{291} & \textbf{1948} &     &     &    &    &  \\
13 & \textbf{0} & \textbf{1} & \textbf{9}  & \textbf{246} & \gr{1300} &     &     &    &    &  \\
14 & \textbf{0} & \textbf{1} & \textbf{13} & \textbf{340} & \gr{1323} &     &     &    &    &  \\
15 & \textbf{0} & \textbf{1} & \textbf{16} & \textbf{349} & \gr{1550} &     &     &    &    &  \\
16 & \textbf{0} & \textbf{1} & \textbf{25} & \textbf{525} & \gr{1925} &     &     &    &    &  \\
\hline
\end{tabular}
\caption{Number of fusion rings available on the AnyonWiki per rank (column) and multiplicity (row). A number is bold if all rings with that rank and multiplicity are available. A number is gray if there might be more fusion rings with the specific rank and multiplicity. No number is displayed if we have no fusion rings with that rank and multiplicity.}\label{tab:numpermultrank}
\end{table}

Every fusion ring has a dedicated page where visitors can find information on important properties such as,
\begin{itemize}
  \item A variety of \textbf{common names} of the fusion ring, together with its \textbf{AnyonWiki code} of four integers and \textbf{UUID}. 
  \item The \textbf{multiplication table} of the fusion ring. This table implicitly settles conventions for the naming and order of the basis elements of the ring. 
  \item The \textbf{Frobenius--Perron dimensions} of the basis elements of the ring and the Frobenius--Perron dimension of the ring itself.
      \item A variety of \textbf{boolean properties} that denote whether the ring is commutative, (weakly) integral, simple, nilpotent, a group ring, etc.
  \item The \textbf{automorphism group} of the fusion ring together with a minimal size set of generators. These are expressed as permutations of the basis elements of the fusion ring using cycle notation. 
  \item The non-trivial {fusion subrings} of the fusion ring.
  \item All possible {gradings} of the fusion ring by a finite group.
  \item The \textbf{characters} of the fusion ring, if it is commutative.
  \item The \textbf{formal codegrees} of the fusion ring, if it is commutative.
  \item The \textbf{tensor product decompositions} of the fusion ring into fusion subrings.
\end{itemize}

Just as for the individual pages on fusion categories, the data featured for fusion rings can be expanded in the future.

\section{The AnyonWikiDatabase}
\label{sec: anyonwikidatabase}
The data that is displayed on the AnyonWiki is automatically fetched from the open-source AnyonWikiDatabase Github repository, which is backed up by the Software Heritage Foundation.\footnote{The current revision of the AnyonWikiDatabase is located at the following URL \url{https://archive.softwareheritage.org/swh:1:rev:41c1452cf6d24d2ea3b8f5819d974f0afd2a1d50;origin=https://github.com/anyonwiki/AnyonWikiDatabase;visit=swh:1:snp:eaa28f68495d508282146e256a30e8e629b14b79}} At the moment this repository provides the data on the fusion rings and fusion categories to AnyonWiki.org. It also contains the algebraic numbers, together with various  \LaTeX\ representations, necessary for expressing the properties of the fusion rings and categories. As more data is generated or contributed to the project the AnyonWikiDatabase will be updated accordingly. The data on the fusion rings (resp. categories) is stored as standard JSON files~\cite{ecma404} containing the following three fields
\begin{itemize}
\item \texttt{data}: A JSON dictionary mapping unique UUIDs of fusion rings (resp. categories) to JSON dictionaries containing the data of each fusion ring (resp. category).
\item \texttt{info}: A full description of the formatting conventions used for each piece of data stored as part of a fusion ring (resp. category). As the database gets updated and fusion rings (resp. categories) get more properties, this description also gets updated so that it is always clear how to interpret each piece of data.
\item \texttt{order}: A list of UUIDs of fusion rings (resp. categories) in the order they were exported. This field is meant to be used for testing purposes for programming languages that import JSON files as unordered dictionaries.
\end{itemize}

The design choices behind how data is formatted in the AnyonWikiDatabase follows the relevant parts of the FAIR (Findable, Accessible, Interoperable, Reusable) data principles~\cite{wilkinson2016fair}. In particular:
\begin{itemize}
\item The fusion rings, categories, and any other supporting data and software use git as their version control system. Therefore, any changes can be easily tracked and, by using commit messages, we can inform users when errors have been fixed.
\item The AnyonWikiDatabase is backed up by the Software Heritage Foundation~\cite{abramatic2018building,di2017software} so it remains accessible as long as possible. 
\item Instructions for how to interpret the data are an intrinsic part of the data itself. This way we avoid issues where the data gets separated from the interpretation and researchers have to guess how to interpret it. 
\item The data of the fusion rings and categories is expressed using standard data types that any programming language can import. For example, complex floating-point numbers are stored as arrays of two real floating points numbers, $[a,b]$, from which the complex number $a+ib$ can be reconstructed. Algebraic numbers are stored as strings of the form \texttt{"a0\_a1\_...\_an\_\_r"} where \texttt{"a0"},$\ldots$,\texttt{"an"} are the coefficients of the polynomial $a_0 + a_{1}x + \cdots + a_{n} x^{n}$ and $r$ is the index of its root according to a well-defined ordering which is explained as part of the \texttt{info} field of the JSON file. This way, any mathematical software package can import our data numerically, and, if supported, symbolically as well. 
\item The data of the AnyonWiki is findable. Every fusion ring and fusion category is identified by a UUID which stored as part of its own entry. UUIDs of fusion rings and categories can be used to find them on the AnyonWiki and each fusion ring and fusion category has a dedicated page which also lists their UUID. Because snapshots of the full database are archived by Software Heritage, they can be referred to by persistent identifiers (SWHIDs). Therefore, any version of any piece of data can be pointed to unambiguously.
\end{itemize}

\begin{remark}{Why the AnyonWiki uses JSON files and not \texttt{.mrdi} files}{}
  The \texttt{.mrdi} file format~\cite{DVJL24,MaRDI22} is a JSON-based format for storing algebraic structures. It is the default serialization format used by OSCAR~\cite{OSCAR-book,JKL25,DV26} and is already used to publish research data, see e.g.~\cite{GZ25}. Every \texttt{.mrdi} file records the software system and version that wrote it. Interestingly, the format is designed such that the user can reconstruct the full algebraic context of the stored objects, which is useful when you have, e.g., a chain of algebraic constructions that depend on each other such as a custom matrix algebra over a custom ring over a custom number field. Since our data only involves lists of algebraic numbers, whose construction is straightforward, we opted instead for a representation that is independent of any particular software package. While this representation prioritizes readability and interoperability over import speed, it is easier to access and researchers are free to convert the data of the AnyonWikiDatabase into any format they like, including \texttt{.mrdi}.
\end{remark}

\section{Citing the AnyonWiki and its contributing sources}
\label{sec: citing}
If you use the AnyonWiki in your research, even if only to get an overview of fusion categories with certain properties, to quickly find out some numeric property about a fusion ring, or to check whether a fusion ring is categorifiable or not, we kindly ask you to cite this paper.

Citing \emph{only} this paper alone, however, is often not sufficient since the data of the AnyonWiki comes from a variety of publications in combination with respective software packages. Since multiple authors have contributed varying amounts of time and effort, the AnyonWiki and the AnyonWikiDatabase provide a system that allows one to cite the relevant sources in a fair manner.

In particular, the AnyonWikiDatabase stores specific references to individual pieces of data via the fields \texttt{references} and \texttt{software}.
The \texttt{references} field contains a dictionary that maps names of data fields of the fusion ring (resp. category) to the references used to obtain that data. 
The \texttt{software} field does the same but maps to the relevant revisions of software packages used to obtain the data.

It is important to note that the AnyonWiki only provides references that played a significant role in providing the data in the way it is stored in the AnyonWikiDatabase. In particular, we do not include references to papers on theoretical properties of fusion rings and categories, nor do we provide references to papers that provide data that did not lead to the specific data stored in the AnyonWikiDatabase. This is to avoid having to add hundreds of references to, e.g., the page on the Fibonacci or Ising fusion ring.

For example, at the moment, the JSON file for the Ising fusion ring contains the following entries 
\begin{verbatim}
"references": {
  "all": ["https://doi.org/10.1063/5.0148848"]
},
"software": {
  "all": {
    "doi": null,
    "name": "FusionRings",
    "swhid": "swh:1:rev:0dc2a749f029b30ca0e4c132d63432b4031ca68b;
              origin=https://github.com/anyonwiki/FusionRings;
              visit=swh:1:snp:efd4aefa347ea7a7b91d0de2da9a41ad47ffea09",
    "url": "https://github.com/anyonwiki/FusionRings",
    "version": "0.2.5"
  }             
}
\end{verbatim}
where we split the software heritage ID (SWHID) in multiple lines to make it more readable. The references field points to the paper(s) that introduced the software that contributed to finding the initial data on the fusion rings. The software field points to a JSON dictionary that specifies the specific revision of the software that generated the data. Any listed piece of software must have either a DOI or a SWHID that provides a permanent link to the specific version of the software used to obtain the data. This is to protect the database and website from containing broken links. 
If extra data from different sources is added, then (a) the new data would get a field name, say \texttt{newprop}, the papers that led to the software would be added to the references field via   \texttt{"newprop" : [ref\_1,...,ref\_n]}, (b) the field \texttt{all} would be changed to \texttt{all\_other\_data}, and (c) a similar change would be made to the \texttt{software} field. Since software changes frequently, it is often the case that the software field has more, and more specific, references than the \texttt{references} field. For example, say that a new field, \texttt{newprop} were added to a fusion ring and
\begin{itemize}
\item the paper that provided the algorithm (or software) to compute \texttt{newprop} has DOI \texttt{https://doi.org/a/b}, and
  \item the software that provided the data is \texttt{newsoft} version \texttt{1.0} with no SWHID, DOI \texttt{https://somesite/c/d}, and url \texttt{https://newsoft.org}, 
\end{itemize}
then that specific fusion ring would have the following software and references fields 
\begin{verbatim}
"references": {
  "all_other_data": ["https://doi.org/10.1063/5.0148848"],
  "newprop": ["https://doi.org/a/b"]
},
"software": {
  "all_other_data": {
    "doi": null,
    "name": "FusionRings",
    "swhid": "swh:1:rev:0dc2a749f029b30ca0e4c132d63432b4031ca68b;
                origin=https://github.com/anyonwiki/FusionRings;
                visit=swh:1:snp:efd4aefa347ea7a7b91d0de2da9a41ad47ffea09",
    "url": "https://github.com/anyonwiki/FusionRings",
    "version": "0.2.5"
  },
  "newprop": {
    "doi": "https://somesite/c/d",
    "name": "newsoft",
    "swhid": null,         
    "url": "https://newsoft.org",
    "version": "1.0" 
  }             
}
\end{verbatim}

These references are automatically fetched by AnyonWiki.org and each individual page on a fusion ring (or category) displays the information under the ``References" and ``Software" sections. By attaching references to each individual piece of data we reduce the risk that the AnyonWiki becomes a black hole of citations, i.e., that researchers only cite the AnyonWiki but not the papers that lead to the data they used.

The AnyonWikiDatabase itself (the GitHub repository) should not be cited since it contains contributions from various sources but is only maintained by several people. Instead (a) if you need to cite the database as a whole, cite this paper, or (b) if you use specific data, cite the respective entries in its \texttt{references} and \texttt{software} fields in addition to this paper.

\section{Future directions and contributing to AnyonWiki.org}
\label{sec: future and contributing}
The AnyonWiki, the AnyonWikiDatabase, and their supporting packages are living and breathing scientific projects that are continuously upgraded and expanded. We will add more data to the AnyonWikiDatabase such as, for example, the centers of all the  multiplicity-free fusion categories up to rank $5$, as  computed in~\cite{center-computation}. We will also continue developing software to compute and work with data on fusion categories. At the moment, a successor to Anyonica~\cite{anyonica} --- the package that provided the $F$-, $P$-, and $R$-symbols of the multiplicity-free fusion categories up to rank $7$ --- is under development. The new package, Lyctor, is an open-source Julia package that will allow researchers to work with skeletal fusion categories. Lyctor will give direct access to the fusion categories stored in the AnyonWikiDatabase and provide functions to compute and probe their properties as well as to categorify fusion rings.
Besides software and data we also plan to continue upgrading the user interface to provide more powerful tools to the end-user without sacrificing user-friendliness. 

\subsection*{Contributing to AnyonWiki.org}
Anyone is welcome to contribute data to the AnyonWiki. Proposing to add data can be done by contacting the AnyonWiki team via the contact form on AnyonWiki.org, by emailing \href{mailto:contact@anyonwiki.org}{contact@anyonwiki.org}, or by contacting one of the authors of this paper directly. 

However, there are certain minimal requirements for data to be included. Any piece of data should come either (a) with a DOI referring to a paper or (b) with a dictionary containing information on the software used to generate it. For the latter we require at least a DOI or a SWHID, the software name, version, and website. 

Data for fusion rings and fusion categories have specific extra requirements. The minimum amount of information required for each structure is as follows.
\begin{itemize}
\item For a fusion ring, we require the fusion coefficients, i.e.~the structure constants of the fusion ring.
\item For a fusion category, we require the $F$-symbols, i.e.~the solutions to the pentagon equations. At this time the project only accepts symbolic algebraic (along with numerical) $F$-symbols, as opposed to numerical $F$-symbols alone. 
\end{itemize}
The data should always be formatted according to the conventions set out in the \texttt{info} field of the JSON file that contains the current list of such structures.

\printbibliography
\end{document}